\documentclass[12pt]{article}

\usepackage{amsfonts}
\usepackage{amsmath}
\usepackage{amsthm}

\usepackage{tikz}
\usetikzlibrary{calc}

\newtheorem{fact}{Fact}

\usepackage[algo2e,ruled,vlined]{algorithm2e}

\theoremstyle{definition}
\newtheorem{example}{Example}

\DeclareMathOperator{\rank}{rank}
\DeclareMathOperator{\mr}{mr}

\usepackage{graphicx}
\usepackage{amssymb}
\usepackage{epstopdf}
\usepackage{tikz-feynman}
\usepackage{hyperref}
\usepackage{graphicx} 
\usepackage{color}
\usepackage{tikz}
\usepackage{mathrsfs}
\usepackage{float}
\usepackage{color}
\usepackage{pagecolor}
\usepackage{subcaption}

\definecolor{grayy}{HTML}{888888}

\title{An exact algorithm for the minimum rank of a graph}

\author{Boris Brimkov\thanks{Department of Mathematics and Statistics, Slippery Rock University, Slippery Rock, PA  (boris.brimkov@sru.edu).} \and 
Zachary Scherr \thanks{Department of Mathematical Sciences, Susquehanna University,  Selinsgrove, PA (scherr@susqu.edu)}
}

\renewenvironment{thebibliography}[1]{
  \begin{oldthebibliography}{#1}
    \setlength{\itemsep}{0em}
    \setlength{\parskip}{0em}
}
{
  \end{oldthebibliography}
}

\date{}

\begin{document}
\maketitle

\abstract{The minimum rank of a graph $G$ is the minimum rank over all real symmetric matrices whose off-diagonal sparsity pattern is the same as that of the adjacency matrix of $G$. In this note we present the first exact algorithm for the minimum rank of an arbitrary graph $G$. In particular, we use the notion of determinantal rank to transform the minimum rank problem into a system of polynomial equations that can be solved by computational tools from algebraic geometry and commutative algebra. We provide computational results, explore possibilities for improvement, and discuss how the algorithm can be extended to other problems such as finding the minimum positive semidefinite rank of a graph.
}

\section{Introduction}

Let $S_n(\mathbb{R})$ denote the set of real symmetric $n\times n$ matrices. For a matrix $A\in S_n(\mathbb{R})$, $\mathcal{G}(A)$ denotes the graph with vertex set $\{1,\ldots,n\}$ and edge set $\{\{i,j\}:A_{ij}\neq 0,\; 1\leq i<j\leq n\}$. Note that the diagonal of $A$ is not used when constructing $\mathcal{G}(A)$. The set of \emph{symmetric matrices associated with a graph} $G$ is defined as $\mathcal{S}(G)=\{A\in S_n(\mathbb{R}):\mathcal{G}(A)=G\}$. The \emph{minimum rank} of $G$ is defined as $\mr(G)=\min\{\text{rank}(A):A\in\mathcal{S}(G)\}$. The minimum rank problem is a special case of the matrix completion problem which has numerous theoretical and practical applications (such as the million-dollar Netflix challenge~\cite{netflix}); it is also related to the inverse eigenvalue problem \cite{hogben3}, quantum controllability on graphs \cite{godsil},  and various other problems in spectral graph theory and combinatorial matrix theory.

The minimum rank problem was first studied in 1996 by Nylen \cite{nylen}, who gave an algorithm for computing the minimum rank of trees; this algorithm was later improved in \cite{johnson1,johnson2,wei}, and generalized to block-cycle graphs in \cite{barioli1}. The graphs having very large and very small minimum ranks have been characterized in \cite{barrett2,barrett,barrett1,hogben1,johnson3}. Decomposition formulas have been derived for computing the minimum ranks of graphs with cut vertices \cite{barioli2,Hsieh} and joins of graphs \cite{barioli3} in terms of the minimum ranks of certain subgraphs. The effects of edge subdivisions \cite{barrett3,barrett4}, edge deletions \cite{Edholm}, and graph complements \cite{barioli4,hogben2} on the minimum rank have also been explored. Upper and lower bounds for the minimum rank of a graph can be obtained using graph theoretic parameters such as the zero forcing number and its variants \cite{AIM-workshop,barioli7,brimkov,brimkov2,gentner,huang}, Colin de Verdi\`{e}re type parameters \cite{barioli7,barioli5,hogben1}, ordered and induced subgraphs \cite{mitchell}, and other methods. Techniques for computing the minimum rank of small graphs are described in \cite{DeLoss1}, and are combined in \cite{DeLoss2} with the bounds mentioned above to compute the minimum ranks of all graphs on up to 7 vertices. See \cite{barioli_zf_mr,survey} for a survey of recent results on the minimum rank problem.

Despite this extensive research, the literature notably lacks an exact algorithm for computing the minimum rank of an arbitrary graph in finite time; researchers in the field have said it would be ``incredibly valuable if such a thing exists" \cite{hogben_letter}. In this note, we present such an algorithm using two well-known facts from linear algebra and commutative algebra. 
We demonstrate the capabilities and limitations of our algorithm by computing the minimum ranks of several graphs, some of which could previously not be computed by any automated method. Some possibilities for improvement, extensions to other problems, and directions for future work are also discussed.

\section{Main results}
\label{section_main}
Given an $n\times m$ matrix $A$, a \emph{minor} of $A$ with \emph{order} $k$ (or $k$-\emph{minor}, for short) is the determinant of a $k\times k$ submatrix obtained from $A$ by removing some of its rows and columns. A \emph{system of polynomial equations} is a set of simultaneous equations $\{p_1(\vec{x})=\ldots =p_m(\vec{x})=0\}$ where $p_i(\vec{x})$, $1\leq i\leq m$, is a polynomial with rational coefficients in several variables $\vec{x}=[x_1,\ldots,x_n]$. A \emph{solution} of a system of polynomial equations is a set of values for $\vec{x}$ which make all equations true.

The first idea we use in our algorithm is the well-known (yet often forgotten) \emph{determinantal rank} of a matrix, i.e., the order of its largest non-vanishing minor. More precisely, we use the fact that the rank of a matrix $A$ is equal to the largest order of any non-zero minor of $A$:

\begin{fact}
\label{fact1}
For any matrix $A$, $\rank(A)= r$ if and only if all $(r+1)$-minors of $A$ are 0 and not all $r$-minors of $A$ are 0.
\end{fact}
\noindent The second idea we use concerns the solution of a system of polynomial equations. It follows from the Tarski-Seidenberg theorem \cite{tarski} that the problem of determining whether a system of polynomial equations has a real solution is decidable:
\begin{fact}
\label{fact2}
It can be determined in finite time whether a system of polynomial equations $\{p_1(\vec{x})=\ldots =p_m(\vec{x})=0\}$ has a real solution. 
\end{fact}
\noindent There are several \emph{tour de force} algorithms for finding the solution of a system of polynomial equations in finite time, if one exists. These include Collins' algorithm for cylindrical algebraic decomposition \cite{collins}, Buchberger's algorithm for computing Gr\"{o}bner bases \cite{buchberger}, and the critical points method of Grigorev and Vorobjov \cite{grigorev}. See \cite{basu} for a detailed survey of algorithms and complexity results on solving systems of polynomial equations. Such algorithms are implemented in computer algebra systems like \emph{Mathematica}, \emph{Sage}, \emph{Maple}, and \emph{Magma}, and dedicated solvers like \emph{Bertini} and \emph{PHCpack}.

We are now ready to describe our algorithm for computing the minimum rank of a graph $G$. Let $A^*$ be a matrix that achieves $\mr(G)$. Since $A^*\in \mathcal{S}(G)$, we can represent the diagonal entries of $A^*$ with variables $\vec{x}=[x_1,\ldots,x_n]$, and the off-diagonal  non-zero entries in the upper triangle of $A^*$ with variables $\vec{y}=[y_1,\ldots,y_t]$. Since $A^*$ is symmetric, the off-diagonal non-zero entries in the lower triangle of $A^*$ can also be represented by $\vec{y}$. Then, a minor of $A^*$ is simply a polynomial with integer coefficients in the variables $\vec{x}, \vec{y}$. Let $\{f_1(\vec{x},\vec{y}),\ldots,f_p(\vec{x},\vec{y})\}$ be the set of all $k$-minors of $A^*$ for some $k\geq 1$. 

By Fact \ref{fact1}, if the system of polynomials $\{f_1(\vec{x},\vec{y})=\ldots=f_p(\vec{x},\vec{y})=0\}$ has a real solution in which $y_1\neq 0,\ldots,y_t\neq 0$, then $\mr(G)\leq k-1$. 
Note that for $1\leq i\leq t$, $y_i\neq 0$ if and only if there exists a $\hat{y}_i$ such that $y_i\hat{y}_i=1$. Thus, $\{f_1(\vec{x},\vec{y})=\ldots=f_p(\vec{x},\vec{y})=0\}$ has a real solution in which $y_1\neq 0,\ldots,y_t\neq 0$ if and only if 
\begin{equation}
\label{reqs}
\{f_1(\vec{x},\vec{y})=\ldots=f_p(\vec{x},\vec{y})=y_1\hat{y}_1-1=\ldots=y_t\hat{y}_t-1=0\}
\end{equation}
has a real solution. By Fact \ref{fact2}, it can be determined whether the system \eqref{reqs} has a real solution in finite time. If this procedure is repeated for $k\geq 1$ and terminated as soon as a real solution to the corresponding system of polynomial equations is found, then the iteration at which a real solution is found is equal to $1+\rank(A^*)=1+\mr(G)$. Note that this procedure always terminates, since $\mr(G)\leq n$. We formally summarize this argument in Algorithm \ref{alg1}.

\vspace{7pt}

\LinesNumbered
\begin{algorithm2e}[h]
\textbf{Input:} A graph $G$ of order $n$\;
\textbf{Output:} $\mr(G)$\;
$\vec{x}=[x_1,\ldots,x_n]\leftarrow$ the diagonal entries of a matrix $A^*\in\mathcal{S}(G)$\;
$\vec{y}=[y_1,\ldots,y_t]\leftarrow$ the non-zero off-diagonal entries in the upper triangle of $A^*$\;
\For{$k=1,\ldots,n$}{
$\{f_1(\vec{x},\vec{y}),\ldots,f_p(\vec{x},\vec{y})\}\leftarrow$ the set of $k$-minors of $A^*$\;
\If{ $\{f_1(\vec{x},\vec{y})=\ldots=f_p(\vec{x},\vec{y})=y_1\hat{y}_1-1=\ldots=y_t\hat{y}_t-1=0\}$ \emph{has a solution}}{
\Return{$\mr(G)\leftarrow k-1$\;}
}
}
\caption{Minimum rank of a graph}
\label{alg1}
\end{algorithm2e}

\vspace{29pt}

\begin{example}
\label{example1}
We illustrate Algorithm \ref{alg1} by applying it to a familiar graph. Let $G$ be the path $P_4$ and let $A^*\in \mathcal{S}(G)$ be a matrix such that $\rank(A^*)=\mr(G)$. Then,

\[ A^*=\left( \begin{array}{cccc}
a & b & 0 &0 \\
b & c & d &0 \\
0 & d & e &f \\
0 & 0 & f &g\end{array} \right),\qquad \text{where}\; b,d,f\neq 0.\] 
Suppose we are in the third iteration of Algorithm \ref{alg1}, i.e., that we have already examined all 1-minors and 2-minors of $A^*$ and found that there is no assignment of variables that yields a rank 2 matrix. Next we consider the 3-minors of $A^*$, which are 
\begin{equation*}
\left| \begin{array}{ccc}
c & d & 0 \\
d & e & f \\
0 & f & g \end{array} \right|, 
\left| \begin{array}{ccc}
b & d & 0 \\
0 & e & f \\
0 & f & g \end{array} \right|, 
\left| \begin{array}{ccc}
b & c & 0 \\
0 & d & f \\
0 & 0 & g \end{array} \right|, 
\left| \begin{array}{ccc}
b & c & d \\
0 & d & e \\
0 & 0 & f \end{array} \right|, 
\end{equation*}
\begin{equation*}
\left| \begin{array}{ccc}
b & 0 & 0 \\
d & e & f \\
0 & f & g \end{array} \right|,
\left| \begin{array}{ccc}
a & 0 & 0 \\
0 & e & f \\
0 & f & g \end{array} \right|,
\left| \begin{array}{ccc}
a & b & 0 \\
0 & d & f \\
0 & 0 & g \end{array} \right|,
\left| \begin{array}{ccc}
a & b & 0 \\
0 & d & e \\
0 & 0 & f \end{array} \right|,
\end{equation*}
\begin{equation*}
\left| \begin{array}{ccc}
b & 0 & 0 \\
c & d & 0 \\
0 & f & g \end{array} \right|,
\left| \begin{array}{ccc}
a & 0 & 0 \\
b & d & 0 \\
0 & f & g \end{array} \right|,
\left| \begin{array}{ccc}
a & b & 0 \\
b & c & 0 \\
0 & 0 & g \end{array} \right|,
\left| \begin{array}{ccc}
a & b & 0 \\
b & c & d \\
0 & 0 & f \end{array} \right|,
\end{equation*}
\begin{equation*}
\left| \begin{array}{ccc}
b & 0 & 0 \\
c & d & 0 \\
d & e & f \end{array} \right|,
\left| \begin{array}{ccc}
a & 0 & 0 \\
b & d & 0 \\
0 & e & f \end{array} \right|,
\left| \begin{array}{ccc}
a & b & 0 \\
b & c & 0 \\
0 & d & f \end{array} \right|,
\left| \begin{array}{ccc}
a & b & 0 \\
b & c & d \\
0 & d & e \end{array} \right|,
\end{equation*}
where the minor in the $i^\text{th}$ row and $j^\text{th}$ column in the array above is obtained by deleting the $i^\text{th}$ row and $j^\text{th}$ column of $A^*$. Evaluating the determinants, we obtain
\begin{equation*}
\begin{array}{cccc}
ceg-cf^2-d^2g, & beg-bf^2, & bdg, & bdf,\\
beg-bf^2, & aeg-af^2, & adg, & adf,\\
bdg, & adg, & acg-b^2g, & acf-b^2f,\\
bdf, & adf, & acf-b^2f, & ace-ad^2-b^2e.
\end{array}
\end{equation*}
We set each of these determinants equal to zero. Moreover, to ensure that the variables $b$, $d$, and $f$ are nonzero, we introduce three new variables $\hat{b}$, $\hat{d}$, and $\hat{f}$ and  the equations $b\hat{b}=1$, $d\hat{d}=1$, and $f\hat{f}=1$. Finally, we check whether the resulting system of polynomial equations has a solution:
\begin{eqnarray*}
ceg-cf^2-d^2g=beg-bf^2=bdg=bdf=beg-bf^2=aeg-af^2=&\\
adg=adf=bdg=adg=acg-b^2g=acf-b^2f=bdf=adf=&\\
acf-b^2f=ace-ad^2-b^2e=b\hat{b}-1=d\hat{d}-1=f\hat{f}-1=&0.
\end{eqnarray*}
The Gr\"{o}bner basis of these polynomials (computed by the computer algebra system \emph{Mathematica}) contains 1, which means the system of polynomial equations has no solution. Hence, there is no assignment of the variables $a,b,c,d,e,f,g$ in $A^*$ which produces a matrix of rank 2. Repeating this procedure with the 4-minor of $A^*$ yields the system of polynomials 
\begin{equation*}
b^2 f^2 - a c f^2 - a d^2 g - b^2 e g + a c e g=b \hat{b}-1 = d \hat{d}-1  =f \hat{f} -1 = 0.
\end{equation*}
The Gr\"{o}bner basis of these polynomials (again computed by \emph{Mathematica}) does not contain 1, and the two families of solutions are
\begin{equation*}
\left\{e=\frac{b^2 f^2 - a c f^2 - a d^2 g}{(b^2 - a c) g}, \hat{b}=\frac{1}{b}, \hat{d}= \frac{1}{d}, \hat{f}= \frac{1}{f}\right\},
\end{equation*}
\begin{equation*}
\left\{c =\frac{b^2}{a}, g =0, \hat{b}=\frac{1}{b}, \hat{d}= \frac{1}{d}, \hat{f}= \frac{1}{f}\right\}.
\end{equation*}
Then, using the first family of solutions and choosing
$a=c=0$, $b=d=e=f=g=1$, we obtain the following rank 3 matrix that achieves the minimum rank of $G$:
\begin{align*}
\begin{pmatrix}
0 & 1 & 0 &0 \\
1 & 0 & 1 &0 \\
0 & 1 & 1 &1 \\
0 & 0 & 1 &1
\end{pmatrix}.
\end{align*}
\end{example}

\subsection{Computational results}

We now illustrate the 
scope of Algorithm \ref{alg1} by applying it to some less familiar graphs. As a baseline, we consider the minimum rank program developed by DeLoss et al. in \cite{DeLoss1,DeLoss3,DeLoss2}; this program computes several combinatorial upper and lower bounds for $\mr(G)$ and returns $\mr(G)$ if some lower bound equals some upper bound. The program also attempts connected component and cut vertex decompositions, leveraging the fact that $\mr(G)$ can be expressed in terms of the minimum ranks of the connected and biconnected components of $G$ (or related subgraphs; see \cite{barioli2} for more details). Finally, if $G$ is a tree, the program uses an exact algorithm for trees to compute $\mr(G)$ (cf. \cite{AIM-workshop}). If none of the above methods succeed in computing $\mr(G)$, then the program returns the best upper and lower bounds for $\mr(G)$. 

We tested Algorithm \ref{alg1} on graphs whose minimum ranks could not be found by the program of DeLoss et al. The minimum ranks of these graphs are known, but were obtained by manually finding matrices in $\mathcal{S}(G)$ whose rank matched a lower bound for $\mr(G)$, or by inspecting the structure of $G$ and using combinatorial arguments on a case-by-case basis to deduce $\mr(G)$ (see Proposition 4.1 in \cite{DeLoss1} and Table 2 in \cite{DeLoss2}). In contrast, Algorithm \ref{alg1} computes the minimum ranks of these graphs without any human intervention. The results and runtimes are reported in Table \ref{table1}. The names of the graphs are their Atlas numbers and their adjacencies can be found in \emph{The Atlas of Graphs} \cite{read}. The computations were performed on a Lenovo Thinkpad with a 2.80GHz Intel i7-7700HQ CPU and 8GB of RAM, running Mathematica 10.0 on Windows 10.

\begin{table}

\def\arraystretch{1.1}
\begin{center}
\caption{Minimum ranks of some graphs}
\begin{tabular}{ccccr}
\hline
$G$ & $|V|$ & $|E|$ & $\mr(G)$ & time (s)\\
\hline
558	&	7	&	9	&	3	&	37.4	\\
669	&	7	&	10	&	3	&	198.8	\\
678	&	7	&	10	&	3	&	267.3	\\
679	&	7	&	10	&	4	&	40310.4	\\
721	&	7	&	10	&	3	&	82.6	\\
791	&	7	&	11	&	3	&	2202.2	\\
801	&	7	&	11	&	3	&	2831.6	\\
812	&	7	&	11	&	3	&	155.2	\\
831	&	7	&	11	&	3	&	392.3	\\
832	&	7	&	11	&	3	&	644.3	\\
846	&	7	&	11	&	3	&	957.8	\\
\hline
\end{tabular}\hspace{30pt}
\begin{tabular}{ccccr}
\hline
$G$ & $|V|$ & $|E|$ & $\mr(G)$ & time (s)\\
\hline
863	&	7	&	11	&	3	&	4682.5	\\
873	&	7	&	11	&	3	&	4125.9	\\
878	&	7	&	11	&	3	&	3794.7	\\
913	&	7	&	12	&	3	&	30178.7	\\
918	&	7	&	12	&	3	&	2207.3	\\
924	&	7	&	12	&	3	&	5404.9	\\
932	&	7	&	12	&	3	&	6958.3	\\
944	&	7	&	12	&	3	&	11664.0	\\
953	&	7	&	12	&	3	&	40402.5	\\
956	&	7	&	12	&	3	&	36576.5	\\
958	&	7	&	12	&	3	&	32794.7	\\
\hline
\end{tabular}

\label{table1}
\end{center}
\end{table}

As can be seen from Table \ref{table1}, the runtime of Algorithm \ref{alg1} typically increases exponentially with $|E(G)|$. For example, graph 558 has the fewest edges among the graphs tested, and took the least amount of time to be solved; graphs 953, 956, and 958 had three more edges than 558 and took roughly three orders of magnitude longer to be solved. This can be explained by the fact that  polynomial equation solvers typically require $(pd)^{2^{O(v)}}$ time, where $p$ is the number of polynomials in the system, $d$ is the maximum degree of a polynomial in the system, and $v$ is the number of variables (see \cite{collins,wuthrich} and the survey of Ayad \cite{ayad}; singly exponential methods have been proposed, e.g. in \cite{grigorev}, but have some limitations). Thus, if $G$ has more edges, the system of polynomial equations obtained from $G$ in Algorithm~\ref{alg1} has more variables and more equations, and thus takes exponentially longer to solve. The runtime of Algorithm \ref{alg1} also typically increases significantly with the minimum rank of $G$. For example, graph 679 has the same number of vertices and edges as graphs 669, 678, and 721, but has a larger minimum rank, and took two orders of magnitude longer to be solved. This can be explained by the fact that if $G$ has a larger minimum rank, then Algorithm \ref{alg1} has to perform more iterations involving larger $k$-factors, which requires the solution of larger systems of equations with higher degrees. Note that since a symmetric $n\times n$ matrix has $\binom{n}{k}^2$ $k$-minors and at most $\frac{n(n-1)}{2}$ distinct nonzero off-diagonal entries, the $k^{\text{th}}$ iteration of Algorithm \ref{alg1} requires solving a system of $\Omega(n^{2k})$ polynomial equations.

Next, we apply Algorithm \ref{alg1} to two slightly larger graphs -- the paths $P_{11}$ and $P_{12}$. It is well known that $\mr(P_n)=n-1$ for any path $P_n$; nevertheless, $P_{11}$ and $P_{12}$ are illustrative of the algorithm's runtime and output. The runtimes (in seconds) of the iterations of Algorithm \ref{alg1} on $P_{11}$ are as follows:
\[0.0, \; 0.0, \; 0.3, \; 2.7, \; 16.3, \; 95.9, \; 294.3, \; 490.4, \; 248.0, \; 23.9, \; 0.2.\]
As expected, at the $11^{\text{th}}$ iteration, a solution to the corresponding system of polynomials was found. One automatically generated instance of the solution is reported below; it can be verified that the rank of this matrix is indeed 10. 

\begin{equation}
\label{eq1}
\scriptsize
\left(
\begin{array}{ccccccccccc}
 -\frac{22928}{4129} & 4 & 0 & 0 & 0 & 0 & 0 & 0 & 0 & 0 & 0 \\
 4 & -2 & -3 & 0 & 0 & 0 & 0 & 0 & 0 & 0 & 0 \\
 0 & -3 & -1 & 4 & 0 & 0 & 0 & 0 & 0 & 0 & 0 \\
 0 & 0 & 4 & -3 & -4 & 0 & 0 & 0 & 0 & 0 & 0 \\
 0 & 0 & 0 & -4 & 1 & -3 & 0 & 0 & 0 & 0 & 0 \\
 0 & 0 & 0 & 0 & -3 & 0 & -2 & 0 & 0 & 0 & 0 \\
 0 & 0 & 0 & 0 & 0 & -2 & -4 & 2 & 0 & 0 & 0 \\
 0 & 0 & 0 & 0 & 0 & 0 & 2 & 3 & -2 & 0 & 0 \\
 0 & 0 & 0 & 0 & 0 & 0 & 0 & -2 & -3 & -2 & 0 \\
 0 & 0 & 0 & 0 & 0 & 0 & 0 & 0 & -2 & 4 & -2 \\
 0 & 0 & 0 & 0 & 0 & 0 & 0 & 0 & 0 & -2 & -1 \\
\end{array}
\right).
\end{equation}
The runtimes (in seconds) of the iterations of Algorithm \ref{alg1} on $P_{12}$ are as follows:
\[0.0, \; 0.0, \; 1.2, \; 4.8, \; 31.6, \; 207.4, \; 1030.5, \; 2518.6, \; 3072.4, \; 2700.5, \; 150.3, \; 0.6.\]
At the $12^{\text{th}}$ iteration, the following rank 11 matrix realizing $\mr(G)$ was returned: 

\begin{equation}
\label{eq2}
\scriptsize
\left(
\begin{array}{cccccccccccc}
 \frac{114144}{22729} & -3 & 0 & 0 & 0 & 0 & 0 & 0 & 0 & 0 & 0 & 0 \\
 -3 & 0 & -3 & 0 & 0 & 0 & 0 & 0 & 0 & 0 & 0 & 0 \\
 0 & -3 & -4 & -2 & 0 & 0 & 0 & 0 & 0 & 0 & 0 & 0 \\
 0 & 0 & -2 & 3 & 2 & 0 & 0 & 0 & 0 & 0 & 0 & 0 \\
 0 & 0 & 0 & 2 & -3 & -2 & 0 & 0 & 0 & 0 & 0 & 0 \\
 0 & 0 & 0 & 0 & -2 & 2 & -2 & 0 & 0 & 0 & 0 & 0 \\
 0 & 0 & 0 & 0 & 0 & -2 & 1 & 4 & 0 & 0 & 0 & 0 \\
 0 & 0 & 0 & 0 & 0 & 0 & 4 & 3 & -1 & 0 & 0 & 0 \\
 0 & 0 & 0 & 0 & 0 & 0 & 0 & -1 & -2 & -2 & 0 & 0 \\
 0 & 0 & 0 & 0 & 0 & 0 & 0 & 0 & -2 & -1 & 4 & 0 \\
 0 & 0 & 0 & 0 & 0 & 0 & 0 & 0 & 0 & 4 & -3 & -4 \\
 0 & 0 & 0 & 0 & 0 & 0 & 0 & 0 & 0 & 0 & -4 & 1 \\
\end{array}
\right).
\end{equation}
Note that the sequences of runtimes of Algorithm \ref{alg1} for $P_{11}$ and $P_{12}$ are unimodal, and skewed to the left. This can be explained by the fact that the largest number of polynomial equations have to be solved at iteration $\frac{n}{2}$ (since $\arg\max_k\{{n\choose k}^2\}=\frac{n}{2}$), but the maximum degree of the polynomials is achieved at iteration $n$; as noted above, both the number of equations and the maximum degree increase the runtime of the solver.

\section{Discussion and future work}
In this note, we combined two simple and well-known ideas to obtain a long-sought exact algorithm for the minimum rank of a graph. 
Overall, the computational limit of our implementation of Algorithm~\ref{alg1} seems to be at dense graphs with around 7 vertices and sparse graphs with around 9 vertices. For larger graphs, the algorithm runs out of memory or has an impractically long runtime. While the computational scope of Algorithm \ref{alg1} is admittedly limited, it settles the decidability of the minimum rank problem and opens the possibility for improved computational methods and new algorithmic paradigms. For example, in addition to $|E(G)|$ and $\mr(G)$, the runtime of Algorithm~\ref{alg1} depends to a large extent on the efficiency of the subroutine used for solving systems of polynomial equations. While our implementation used a general purpose solver, the systems of polynomial equations that arise from computing determinants of matrices are clearly quite special, and leave much room for specialization and improvement. Thus, an important direction for future work is to develop algorithms that more efficiently determine whether a system of polynomial equations arising from matrix determinants has a solution. Any \emph{a priori} assumptions about the graph whose minimum rank is being computed (e.g., sparsity, symmetry, connectivity) could potentially also be leveraged to obtain a significant speedup. Upper and lower bounds on the minimum rank -- such as the zero forcing number -- remain relevant, since solving systems of polynomial equations typically takes doubly exponential time, whereas computing the zero forcing number and other graph-based bounds on the minimum rank takes (only) exponential time. Thus, the number of iterations of Algorithm~\ref{alg1} can be reduced by having good upper and lower bounds for the minimum rank. Binary search can also be used to further decrease the number of iterations. 

With slight modifications, Algorithm \ref{alg1} can be used to solve several other problems that are related to the minimum rank problem. For example, the \emph{minimum positive semidefinite rank} of $G$, denoted $\mr_+(G)$, is defined as the minimum rank over all positive semidefinite matrices with the same sparsity pattern as $G$. The minimum positive semidefinite rank and similar parameters have  been widely studied, mainly involving characterizations for specific graphs and general bounds (see, e.g.,  \cite{ekstrand,osborne,wang_pos,yang_pos,zimmer}); however, until now, the literature did not contain an exact algorithm for computing these parameters. Given an $n \times n$ matrix $A$ and a set $I\subseteq \{1,\ldots,n\}$ with $|I|=k$, a \emph{principal minor} of $A$ is a minor that corresponds to the rows and columns of $A$ indexed by $I$. By Sylvester's criterion for positive semidefinite matrices \cite{prussing}, a symmetric matrix $A$ is positive semidefinite if and only if all principal minors of $A$ are nonnegative. Thus, augmenting the system of polynomials being solved in line 7 of Algorithm \ref{alg1} with a system of polynomial inequalities dictating that all polynomials corresponding to the principal minors of $A^*$ are nonnegative will ensure that any solution to the system corresponds to a positive semidefinite matrix. Since the Tarski-Seidenberg theorem also applies to systems of polynomial inequalities, Algorithm \ref{alg1} equipped with a solver for polynomial inequalities can be used to compute $\mr_+(G)$. Other minimum rank parameters restricted by the definiteness of the matrix or by other criteria (e.g. zeros on the diagonal, as in  \cite{zero_diag}) can be handled analogously. 

As another example, depending on the polynomial equation solver used, Algorithm~\ref{alg1} could return -- in addition to $\mr(G)$ -- a matrix whose rank equals $\mr(G)$ (as the matrices in \eqref{eq1} and \eqref{eq2} for $P_{11}$ and $P_{12}$), or a characterization of all matrices whose rank equals $\mr(G)$ (as in Example \ref{example1}). This opens possibilities for computationally investigating the set of all matrices that realize the minimum rank of a graph, and finding a matrix that is optimal with respect to certain other criteria.

We also remark that the \emph{maximum nullity} of a graph $G$, defined as $M(G)=\max\{\text{null}(A):A\in\mathcal{S}(G)\}$, and the \emph{maximum multiplicity} of $G$, defined as $mult(G)=\max\{mult_A(\lambda) : A\in \mathcal{S}(G), \lambda\in \mathbb{R}\}$, where $mult_A(\lambda)$ denotes the multiplicity of $\lambda$ as a root of the characteristic polynomial of $A$, are in a sense equivalent to the minimum rank problem and hence can also be solved by Algorithm \ref{alg1}. In particular, since $\text{rank}(A)+\text{null}(A)=n$ for any matrix $A$, it follows that $M(G)=n-\text{mr}(G)$; moreover, since $\lambda$ is an eigenvalue of $A$ if and only if $0$ is an eigenvalue of $A-\lambda I$ and since $mult_0(A)=\text{null}(A)$, it follows that $mult(G)=n-\mr(G)$. Note that the analogously defined maximum rank, minimum nullity, and minimum multiplicity problems are not interesting, since a matrix in $\mathcal{S}(G)$ with rank $n$ can be constructed by choosing each diagonal entry to be greater than the sum of the other entries in its row (by the Gershgorin circle theorem, such a matrix is nonsingular).

Finally, we briefly address the minimum rank problem for other fields. Let $S_n(\mathbb{F})$ denote the set of symmetric $n\times n$ matrices whose entries belong to a field $\mathbb{F}$, let $\mathcal{S}_{\mathbb{F}}(G)=\{A\in S_n(\mathbb{F}):\mathcal{G}(A)=G\}$, and let $\mr_{\mathbb{F}}(G)=\min\{\text{rank}(A):A\in\mathcal{S}_{\mathbb{F}}(G)\}$. For any finite field $\mathbb{F}$, $\mr_{\mathbb{F}}(G)$ can clearly be computed in finite time, since there are a finite number of matrices in $\mathcal{S}_{\mathbb{F}}(G)$. Moreover, variants of Algorithm \ref{alg1} can be used to compute $\mr_{\mathbb{F}}(G)$ for some other infinite fields, such as $\mathbb{C}$; however, it cannot be used to compute $\mr_{\mathbb{Z}}(G)$, since the problem of determining whether a multivariate polynomial equation has an integer solution (Hilbert's tenth problem) is undecidable \cite{matiyasevich}. It would be interesting to determine whether computing $\mr_{\mathbb{Z}}(G)$  is altogether undecidable, or whether there exists a finite time algorithm for it based on a different paradigm. Similarly, it would be interesting to determine whether computing $\mr_{\mathbb{Q}}(G)$  is undecidable; note that the decidability of determining whether a multivariate polynomial equation has a rational solution is still open.

\end{document}